\documentclass[citesort,number,dvips]{arxbj}
\usepackage{dcolumn}
\usepackage{graphicx}

\aid{0}
\volume{17}
\issue{1}
\pubyear{2011}
\firstpage{253}
\lastpage{275}
\doi{10.3150/10-BEJ279}

\makeatletter
\newcolumntype{d}[1]{D{.}{.}{#1}}
 \newtheorem{prop}{Proposition}

\newremark{remark}{Remark}
\newcommand{\iint}{\int\!\!\!\int}
\newcommand{\eqref}[1]{(\ref{#1})}
\renewcommand{\citet}[1]{\citet{#1}}
\renewcommand{\citep}[1]{\citeauthor{#1} \citeyear{#1}}
\makeatother

\begin{document}
\begin{frontmatter}

\title{A goodness-of-fit test for bivariate extreme-value copulas}
\runtitle{Goodness-of-fit testing for extreme-value copulas}

\begin{aug}
\author[a]{\fnms{Christian} \snm{Genest}\thanksref{a}\ead[label=e1]{Christian.Genest@mat.ulaval.ca}\corref{}},
\author[b]{\fnms{Ivan} \snm{Kojadinovic}\thanksref{b}\ead[label=e2]{Ivan.Kojadinovic@univ-pau.fr}},
\author[c]{\fnms{Johanna} \snm{Ne\v{s}lehov\'a}\thanksref{c}\ead[label=e3]{Johanna.Neslehova@math.mcgill.ca}}
\and
\author[d]{\fnms{Jun} \snm{Yan}\thanksref{d}\ead[label=e4]{jun.yan@uconn.edu}}
\runauthor{Genest, Kojadinovic, Ne\v{s}lehov\'a and Yan}
\pdfauthor{Christian Genest, Ivan Kojadinovic, Johanna Neslehov\'a, Jun Yan}
\address[a]{D\'epartement de math\'ematiques et de statistique,
Universit\'e Laval, 1045, avenue de la M\'edecine, Qu\'ebec, Canada G1V 0A6. \printead{e1}}
\address[b]{Laboratoire de math\'{e}matiques et applications, UMR CNRS 5142, Universit\'{e} de
Pau et des Pays de l'Adour, Bo\^{\i}te postale 1155, 64013 Pau cedex, France.
 \printead{e2}}
\address[c]{Department of Mathematics and Statistics, McGill
University, 805, rue Sherbrooke Ouest, Montr\'eal (Qu\'ebec), Canada H3A 2K6. \printead{e3}}
\address[d]{Department of Statistics, University of Connecticut, 215
Glenbrook Road, Storrs, CT 06269, USA. \printead{e4}}
\end{aug}

\received{\smonth{8} \syear{2009}}
\revised{\smonth{1} \syear{2010}}

%
\begin{abstract}
It is often reasonable to assume that the dependence structure of a
bivariate continuous distribution belongs to the class of extreme-value
copulas. The latter are characterized by their Pickands dependence
function. In this paper, a procedure is proposed for testing whether
this function belongs to a given parametric family. The test is based
on a Cram\'er--von Mises statistic measuring the distance between an
estimate of the parametric Pickands dependence function and either one
of two nonparametric estimators thereof studied by Genest and Segers
[\textit{Ann. Statist.} \textbf{37} (2009) 2990--3022]. As the limiting
distribution of the test statistic depends on unknown parameters, it
must be estimated via a parametric bootstrap procedure, the validity of
which is established. Monte Carlo simulations are used to assess the
power of the test and an extension to dependence structures that are
left-tail decreasing in both variables is considered.
\end{abstract}

%

\end{frontmatter}
%

\section{Introduction}
\label{sec:1}

Let $X$ and $Y$ be continuous random variables with cumulative
distribution functions $F$ and $G$, respectively.
Following Sklar \cite{Skl59}, the joint behavior of the pair $(X,Y)$ can be
characterized at every $(x, y) \in\mathbb{R}^2$ by the relation
%
%
\begin{equation}
\label{eq:1}
H(x,y) = \Pr(X \le x , Y \le y ) = C \{ F(x), G(y) \}
\end{equation}
through a unique copula $C$ that captures the dependence between $X$
and $Y$.

When $H$ is known, its marginal distributions can easily be retrieved
from it. The copula can also be readily identified as it is simply the
joint distribution of the pair $(U,V) = (F(X), G(Y))$. In practice,
however, $H$ is often unknown, and the relation between $X$ and $Y$
must be modeled from data.

A copula model for $H$ assumes that equation \eqref{eq:1} holds for
some $F$, $G$ and $C$ from specific parametric classes. This approach
was used, for example, by Frees and Valdez \cite{FreVal98} and Klugman and Parsa \cite{KluPar99} to
analyze data from the Insurance Services Office, Inc. on the indemnity
payment~($X$) and allocated loss adjustment expense ($Y$) for 1500
general liability claims randomly chosen from late settlement lags.
Based on their work and subsequent analysis by other authors, it is
reasonable to assume that for these data, $F$ is inverse paralogistic, $G$ is Pareto
and~$C$ is a Gumbel--Hougaard extreme-value copula.

Extreme-value copulas characterize the limiting dependence structure of
suitably normalized componentwise maxima. They are of special interest
in insurance \cite{DenuitDhaeneGoovaertsKaas2005}, finance \cite{CherubiniLucianoVecchiato2004,McNFreEmb05}
and hydrology \cite{SalDeMKotRos07}, where the occurrence of joint extremes is a risk
management concern.

Pickands \cite{Pic81} showed that if $C$ is a bivariate extreme-value copula, then
%
%
\begin{equation}
\label{eq:2}
C(u,v) = \exp\biggl[ \log(uv) A \biggl\{ \frac{\log(v)} {\log(uv)}
\biggr\} \biggr]
\end{equation}
for all $u, v \in(0,1)$ and a mapping $A\dvtx  [0,1] \to[1/2,1]$, referred
to as the Pickands dependence function, which is convex and such that
$\max(t, 1-t) \leq A(t) \leq1$ for all $t \in[0,1]$. For instance,
an extreme-value copula is said to belong to the Gumbel--Hougaard
family if there exists $\theta\in[1, \infty)$ such that for all $t
\in[0,1$], we have
%
%
\begin{equation}
\label{eq:3}
A(t) = \{ t^\theta+ (1-t)^\theta\} ^{1/\theta}.
\end{equation}

A test that a copula $C$ is of the form (\ref{eq:2}) was developed by
Ghoudi \textit{et al.} \cite{GhoKhoRiv98}; it was recently refined by Ben Ghorbal
\textit{et al.} \cite{BenGenNes09}.
Under the assumption that $C$ is an extreme-value copula, it may be of
interest to check whether its Pickands dependence function $A$ belongs
to a specific parametric class, say $\mathcal{A} = \{A_\theta\dvtx  \theta
\in\mathcal{O} \}$, where $\mathcal{O}$ is an open subset of
$\mathbb
{R}^p$ for some integer $p$.

The purpose of this paper is to examine how the hypothesis $H_0\dvtx  A \in
\mathcal{A}$ can be tested with a random sample $(X_1, Y_1), \ldots,
(X_n, Y_n)$ from $H$. As for all goodness-of-fit tests reviewed by
Berg \cite{Ber09} and Genest \textit{et al.} \cite{GenRemBea09}, the proposed procedure is based
on pseudo-observations $(U_1, V_1), \ldots, (U_n, V_n)$ from copula
$C$, defined, for $i \in\{ 1, \ldots, n \}$, by
%
%
\begin{equation}\label{eq:4}
U_i = F_n (X_i),\qquad V_i = G_n (Y_i),
\end{equation}
where $F_n$ and $G_n$ are rescaled empirical counterparts of $F$ and
$G$, respectively, given by
\[
F_n (x) = \frac{1}{n+1 } \sum_{i=1}^n \mathbf{1} (X_i \leq x),\qquad
G_n (y) =
\frac{1}{n+1} \sum_{i=1}^n \mathbf{1} (Y_i \leq y)
\]
for all $x, y \in\mathbb{R}$. This approach is justified because, as
copulas themselves, the pairs $(U_1, V_1), \ldots, (U_n, V_n)$ of
normalized ranks are invariant under strictly increasing
transformations of $X$ and $Y$. As shown by Kim \textit{et al.} \cite
{KimSilvapulle2007}, it also leads to efficient and robust estimators.

The proposed test is described in Section \ref{sec:2} and its
asymptotic null distribution is given in Section \ref{sec:3}, where a
parametric bootstrap is proposed for the calculation of $P$-values. In
Section \ref{sec:4}, the distributional result is extended to
alternatives that are left-tail decreasing in both variables. This is
instrumental in studying the consistency and power of the test, which
are considered in Sections \ref{sec:5} and \ref{sec:6}, respectively.
The paper concludes with an illustrative example. Technical proofs are
grouped in a series of appendices.

All procedures discussed herein are implemented in the R package
\texttt{copula} \cite{copula} available via the Comprehensive R Archive
Network at \url{http://cran.r-project.org}.

\section{Proposed goodness-of-fit test}
\label{sec:2}

Let $(X_1, Y_1), \ldots, (X_n, Y_n)$ be a random sample from some
unknown continuous bivariate distribution $H$ whose underlying copula
is of the form \eqref{eq:2} with Pickands dependence function $A$. In
order to test the hypothesis
\[
H_0\dvtx  A \in\mathcal{A} = \{ A_\theta\dvtx  \theta\in\mathcal{O} \},
\]
a natural way to proceed is to compare a nonparametric estimator $A_n$
of $A$ to a parametric estimator $A_{\theta_n}$. Several measures of
distance can be used for this purpose, but the Cram\'er--von Mises statistic
%
%
\begin{equation}
\label{eq:5}
S_n = \int_0^1 n | A_n(t) - A_{\theta_n}(t) |^2 \,\mathrm{d}t
\end{equation}
generally leads to more powerful tests than, say, the
Kolmogorov--Smirnov statistic \cite{GenRemBea09}. The choices of
$A_{\theta_n}$ and $A_n$ are discussed next.

\subsection{Parametric estimation of $A$}
\label{sec:2.1}

Under $H_0$, $A_\theta$ may be estimated by $A_{\theta_n}$ using a
consistent estimate $\theta_n$ of $\theta$. Such an estimate can be
derived from the pairs $(U_1, V_1), \ldots, (U_n, V_n)$ via the
maximum pseudo-likelihood method considered by Genest \textit{et al.} \cite{GenGhoRiv95} and
Shih and Louis \cite{ShiLou95}.

To illustrate this approach in a concrete case, let $A_\theta$ be the
generator of the Gumbel--Hougaard copula defined in \eqref{eq:3}. For
all $u, v \in(0,1)$, write
\begin{eqnarray*}
C_\theta(u,v) &= &\exp\biggl[ \log(uv) A_\theta\biggl\{ \frac{\log(v)}
{\log(uv)} \biggr\} \biggr]\\
&=& \exp[ - \{ |\log(u)|^\theta+ |\log
(v)|^\theta\}^{1/\theta} ].
\end{eqnarray*}
As $A_\theta$ is twice differentiable on $(0,1)$, the copula $C_\theta$
has a density given by $c_\theta(u,v) = {\partial^2} C_\theta(u,v) /
{\partial u \,\partial v} $ everywhere on $(0,1)^2$. The maximum
pseudo-likelihood estimator $\theta_n$ is then the value $\theta\in
\mathcal{O} = (1, \infty)$ at which the function
\[
\ell(\theta) = \sum_{i=1}^n \log\{ c_\theta(U_i, V_i)\}
\]
reaches its global maximum. An advantage of this method is that it can
be used even when the parameter space $\mathcal{O}$ is multidimensional.

When $\theta$ is real-valued, a simpler technique which also yields a
consistent estimator is based on the inversion of Kendall's tau. As
shown by Ghoudi \textit{et al.} \cite{GhoKhoRiv98}, the relation
\[
\tau(C) = -1 + 4 \iint_0^{1} C(u,v) \,\mathrm{d}C(u,v) = \int_0^1 \frac
{t(1-t)}{A(t)}\, \mathrm{d}A^\prime(t)
\]
is valid for any extreme-value copula $C$. When $A \in\mathcal{A}$,
$\tau$ is a function of $\theta$ and a rank-based moment estimate of
the latter is obtained by solving the equation $\tau_n = \tau(\theta)$
for $\theta$, where $\tau_n$ is the sample value of Kendall's tau. In
the Gumbel--Hougaard model, for instance, we find $\tau(\theta) =
1-1/\theta$ and hence $\theta_n = \max\{ 1, 1/(1-\tau_n)\}$.

When $\mathcal{O} \subset\mathbb{R}$, we can also obtain consistent,
rank-based estimates of $\theta$ by exploiting its one-to-one
relationship with other nonparametric measures of dependence such as
Spearman's rho, that is,
\[
\rho(C) = -3 + 12 \iint_0^{1} uv \,\mathrm{d}C(u,v) = -1 + \int_0^1 \frac
{1}{ \{
A(t)\} ^2} \,\mathrm{d}t.
\]

\subsection{Nonparametric estimation of $A$}
\label{sec:2.2}

Nonparametric estimators of $A$ are proposed by Genest and Segers \cite{GenSeg09}. For
$i \in\{ 1, \ldots, n \}$, set $\xi_i (0) = -\log(U_i)$, $\xi_i (1)
= -\log(V_i)$ and
\[
\xi_i (t) = \min\biggl\{ \frac{-\log(U_i)}{1-t} , \frac{-\log
(V_i)}{t} \biggr\}
\]
for all $t \in(0,1)$, where $U_i$ and $V_i$ are as in equation \eqref
{eq:4}. Also, let
\[
A_n^{\mathrm{P}}(t) = 1 \Big/ \Biggl\{\frac{1}{n} \sum_{i=1}^n
\xi_i (t)\Biggr\},\qquad
A_n^{\mathrm{CFG}}(t) = \exp\Biggl[ - \gamma- \frac{1}{n}
\sum_{i=1}^n \log
\{ \xi_i(t)\} \Biggr],
\]
where $\gamma= - \int_0^\infty\log(x) \mathrm{e}^{-x} \,\mathrm{d}x
\approx0.577$
is Euler's constant.

The functions $A_n^{\mathrm{P}}$ and $A_n^{\mathrm{CFG}}$ are
rank-based versions of the
estimators of $A$ introduced by Pickands \cite{Pic81} and Cap{\'e}ra{\`a} \textit{et~al.} \cite{CapFouGen97},
respectively. As noted by Genest and Segers \cite{GenSeg09}, these estimators can be
altered to meet the end-point conditions $A_n^{\mathrm{P}}(0) =
A_n^{\mathrm{CFG}}(0) = 1$ and
$A_n^{\mathrm{P}}(1) = A_n^{\mathrm{CFG}}(1) = 1$. However, this
makes no difference asymptotically.

Both $A_n^{\mathrm{P}}$ and $A_n^{\mathrm{CFG}}$ can be expressed as
functionals of the empirical
copula, which may be defined for all $u, v \in[0,1]$ by
\[
C_n(u,v) = \frac{1}{n} \sum_{i=1}^n \mathbf{1}(U_i \leq u, V_i \leq v).
\]
To be specific, the following relations hold for all $t \in[0,1]$:
\begin{eqnarray*}
A_n^{\mathrm{P}}(t) & = & 1 \Big/ \int_0^1 C_n(x^{1-t},x^t)
\frac{\mathrm{d}x}{x}
, \\
A_n^{\mathrm{CFG}}(t)& = & \exp\biggl\{- \gamma+ \int_0^1 \{
C_n(x^{1-t},x^t) -
\mathbf{1}(x > \mathrm{e}^{-1}) \} \frac{\mathrm{d}x}{x \log(x)} \biggr\}
.
\end{eqnarray*}

It was shown by R\"uschendorf \cite{Rue76} that under weak regularity conditions, the
process $\sqrt{n} (C_n - C)$ converges in law to a Gaussian limit
$\mathbb{C}$, that is,\ $\sqrt{n} (C_n - C) \rightsquigarrow
\mathbb{C}$ as $n \to
\infty$. We may thus expect $A_n^{\mathrm{P}}$ and $A_n^{\mathrm
{CFG}}$ to be consistent and
asymptotically Gaussian. This is shown by Genest and Segers \cite{GenSeg09}, provided
that $A$ is twice continuously differentiable. Their Theorem~3.2 states that
\[
\mathbb{A}_n^{\mathrm{P}}= \sqrt{n} (A_n^{\mathrm{P}}- A)
\leadsto\mathbb{A}^{\mathrm{P}},\qquad
\mathbb{A}_n^{\mathrm{CFG}}= \sqrt{n} ( A_n^{\mathrm{CFG}}- A)
\leadsto\mathbb{A}^{\mathrm{CFG}}
\]
as $n \to\infty$ in $\mathcal{C} [0,1]$, where, for all $t \in[0,1]$,
\begin{eqnarray*}
\mathbb{A}^{\mathrm{P}}(t) & = & - A^2(t) \int_0^1 \mathbb
{C}(x^{1-t},x^t) \frac{\mathrm{d}x}{x}
, \\
\mathbb{A}^{\mathrm{CFG}}(t)& = & A(t) \int_0^1 \mathbb
{C}(x^{1-t},x^t) \frac{\mathrm{d}x}{x \log
(x)} .
\end{eqnarray*}
\begin{remark*} Observe that, in principle, the statistics
$S_n^{\mathrm{P}}$ and
$S_n^{\mathrm{CFG}}$ could be extended to arbitrary dimension $d \ge
3$ because
$d$-variate extreme-value copulas are characterized by $(d-1)$-place
Pickands dependence functions \cite{FalRei05}. At present, however,
multivariate analogs of the rank-based estimators $A_n^{\mathrm{P}}$
and $A_n^{\mathrm{CFG}}$ are
unavailable. To see how the estimation can proceed in the $d$-variate
case when the marginal distributions are known, refer to \cite
{ZhaWelPen08} or \cite{GudSeg09}.
\end{remark*}

\section{Asymptotic null distribution of the test statistic}
\label{sec:3}

The asymptotic distribution of the goodness-of-fit statistic $S_n$
depends on the joint behavior of $\Theta_n = \sqrt{n} (\theta_n -
\theta)$ and either $\mathbb{A}_n^{\mathrm{P}}$ or $\mathbb
{A}_n^{\mathrm{CFG}}$ under $H_0$. Suppose that the
class $\mathcal{A} = \{ A_\theta\dvtx  \theta\in\mathcal{O} \}$ satisfies
the following conditions:
\begin{enumerate}
\item[(A)] the parameter space $\mathcal{O} $ is an open subset of
$\mathbb{R}^p$;
\item[(B)] for every $\theta\in\mathcal{O}$, $A_\theta$ is twice
continuously differentiable on $(0,1)$;
\item[(C)] the gradient $\dot A_\theta(t)$ of $A_\theta(t)$ with
respect to $\theta$ satisfies
%
%
\begin{equation}
\label{eq:6}
\lim_{\epsilon\downarrow0} \sup_{\|\theta^* - \theta\| < \epsilon}
\sup_{t \in[0,1]} \|{\dot A}_{\theta^*} (t) - {\dot A}_\theta(t)\|
\rightarrow0,
\end{equation}
where $\| \cdot\| $ denotes the $\ell_2$-norm.
\end{enumerate}

As is proved in Appendix \hyperref[appa]{A}, the process $\mathbb{A}_{n,\theta_n} =
\sqrt{n}
(A_n - A_{\theta_n})$ is then asymptotically Gaussian, both when $A_n =
A_n^{\mathrm{P}}$ and $A_n = A_n^{\mathrm{CFG}}$.

\begin{prop}
\label{prop:1}
Assume $H_0$ holds, that is, $C$ is an extreme-value copula with
Pickands dependence function $A = A_{\theta_0}$ for some $\theta_0
\in
\mathcal{O}$. Further, assume that $\mathcal{A} = \{ A_\theta\dvtx
\theta
\in\mathcal{O} \} $ meets conditions \textup{(A)}--\textup{(C)}.
\begin{enumerate}[(a)]
\item[(a)] If $(\mathbb{A}_n^{\mathrm{P}}, \Theta_n)$
converges to a Gaussian limit
$(\mathbb{A}^{\mathrm{P}},\Theta)$, then $\mathbb{A}_{n,\theta_n}
\leadsto\mathbb{A}^{\mathrm{P}}- {\dot
A}_{\theta_0}^\top\Theta$ as $n \to\infty$ in $\mathcal{C} [0,1]$.
\item[(b)] If $(\mathbb{A}_n^{\mathrm{CFG}}, \Theta_n)$
converges to a Gaussian limit
$(\mathbb{A}^{\mathrm{CFG}},\Theta)$, then $\mathbb{A}_{n,\theta
_n} \leadsto\mathbb{A}^{\mathrm{CFG}}- {\dot
A}_{\theta_0}^\top\Theta$ as $n \to\infty$ in $\mathcal{C} [0,1]$.
\end{enumerate}
\end{prop}

The weak convergence of the statistic $S_n$ defined in \eqref{eq:5}
follows immediately from Proposition~\ref{prop:1} and the continuous
mapping theorem (see, e.g., \cite{vanWel96}, Theorem~1.3.6). As the
limit depends on the unknown parameter value $\theta_0$, we must resort
to resampling techniques to carry out the test. The following
parametric bootstrap procedure can be used to this end. Its validity
depends on regularity conditions adapted from \cite{GenRem08}. These
conditions, listed in Appendix \hyperref[appb]{B}, can be verified for many families of
extreme-value copulas.

\textbf{Parametric bootstrap procedure}

\begin{enumerate}[(1)]
\item[(1)] Compute $A_n$ from the pairs $(U_1, V_1), \ldots, (U_n, V_n)$ of
normalized ranks and estimate $\theta$ using a rank-based estimator, as
discussed in Section \ref{sec:2}.
\item[(2)] Compute the test statistic $S_n$ defined in (\ref{eq:5}).
\item[(3)] For some large integer $N$, repeat the following steps for every
$k \in\{1, \ldots, N\}$:
\begin{enumerate}[(3.1)]
\item[(3.1)] generate a random sample $(X_{1k},Y_{1k}), \ldots,
(X_{nk},Y_{nk})$ from copula $C_{\theta_n}$ and deduce the associated
pairs $(U_{1k},V_{1k}), \ldots, (U_{nk}$, $V_{nk})$ of normalized ranks;
\item[(3.2)] let $A_{nk}$ and $\theta_{nk}$ stand for the versions of
$A_n$ and $\theta_n$ derived from the pairs $(U_{1k},V_{1k}), \ldots,
(U_{nk}, V_{nk})$;
\item[(3.3)] compute
\[
S_{nk} = \int_0^1 n | A_{nk}(t) - A_{\theta_{nk}}(t) |^2 \,\mathrm{d}t.
\]
\end{enumerate}
\item[(4)] An approximate $P$-value for the test is given by $N^{-1} \sum
_{k=1}^N \mathbf{1}(S_{nk} \geq S_n)$.
\end{enumerate}

\section{Extension to left-tail decreasing copulas}
\label{sec:4}

The statistic $S_n$ can be used to build goodness-of-fit tests for the
more general hypothesis
\[
H_0^* \dvtx  C \in\mathcal{C} = \{ C_\theta\dvtx  \theta\in\mathcal{O}\},
\]
where $\mathcal{C}$ is a parametric family of copulas that are
left-tail decreasing (LTD) in both arguments. From \cite{Nelsen2006},
Exercise 5.35, a copula $C$ is LTD in this sense if and only if, for
all $ 0 < u \le u^\prime\le1$ and $0 < v \le v^\prime\le1$,
\begin{equation}
\label{eq:7}
\frac{C(u,v)}{uv} \ge\frac{C(u^\prime, v^\prime)}{u^\prime
v^\prime}
.
\end{equation}
This condition is satisfied for extreme-value copulas, which Garralda-Guillem \cite
{GarraldaGuillem2000} showed to be stochastically increasing in both
variables.

The following result, proved in Appendix \hyperref[appc]{C}, implies that when $C$ is an
LTD copula, $A_n^{\mathrm{P}}$ and $A_n^{\mathrm{CFG}}$ are
consistent, asymptotically Gaussian
estimators of $A_C^\mathrm{P}$ and $A_C^\mathrm{CFG}$, respectively, where,
for all $t \in[0,1]$,
\[
A_C^\mathrm{P} (t) = 1 \Big/ \int_0^1 C(x^{1-t}, x^t) \frac
{\mathrm{d}x}{x}
\]
and
\[
A_C^\mathrm{CFG} (t) = \exp\biggl[ - \gamma+ \int_0^1 \{C(x^{1-t}, x^t) -
\mathbf{1} (x > \mathrm{e}^{-1}) \} \frac{\mathrm{d}x}{x \log(x)} \biggr].
\]

\begin{prop}
\label{prop:2}
Suppose that the copula $C$ has a continuous density and satisfies
condition \eqref{eq:7}. Then $\sqrt{n} (A_n^{\mathrm{P}}-
A_C^\mathrm{P})
\leadsto\mathbb{A}^{\mathrm{P}}_C$ and $\sqrt{n} ( A_n^{\mathrm
{CFG}}- A_C^\mathrm{CFG}) \leadsto
\mathbb{A}^{\mathrm{CFG}}_C$ as $n \to\infty$ in $\mathcal{C}
[0,1]$, where, for all $t
\in[0,1]$,
\begin{eqnarray*}
\mathbb{A}^{\mathrm{P}}_C (t) & = & - \{ A_C^\mathrm{P}(t)\}^2 \int_0^1
\mathbb{C}(x^{1-t},x^t)
\frac{\mathrm{d}x}{x} , \\
\mathbb{A}^{\mathrm{CFG}}_C (t)& = & A_C^\mathrm{CFG} (t) \int_0^1
\mathbb{C}(x^{1-t},x^t) \frac
{\mathrm{d}x}{x \log(x)} .
\end{eqnarray*}
\end{prop}

Incidentally, the mappings $A_C^\mathrm{P}$ and $A_C^\mathrm{CFG}$ are well
defined for any copula $C$, whether or not it is LTD. They reduce to
the Pickands dependence function $A$ when $C$ is of the form \eqref
{eq:2}. Otherwise, they typically differ from one another, but retain
some of the properties of Pickands dependence functions. These facts
are summarized in the following proposition, the proof of which is left
to the reader.

\begin{prop}
\label{prop:3}
Let $C$ be a copula and let $A_C$ denote either $A_C^\mathrm{P}$ or
$A_C^\mathrm{CFG}$. Also, let $W$ and $M$ denote the lower and upper
Fr\'
echet--Hoeffding bounds, respectively. The following statements then hold:
\begin{enumerate}[(a)]
\item[(a)] $A_W (t) \ge A_C (t) \ge A_M(t) = \max(t,1-t)$ for
all $t \in[0,1]$;
\item[(b)] if $C(u,v) \ge uv$ for all $u,v \in[0,1]$, then $A_C
(t) \le1$ for all $t \in[0,1]$;
\item[(c)] if $C(u,v) = C(v,u)$ for all $u,v \in[0,1]$, then
$A_C(t) = A_C(1-t)$ for all $t \in[0,1]$;
\item[(d)] if $C$ is an extreme-value copula with Pickands
dependence function $A$, then $A_C = A$.
\end{enumerate}
\end{prop}

\begin{figure}[b]

\includegraphics{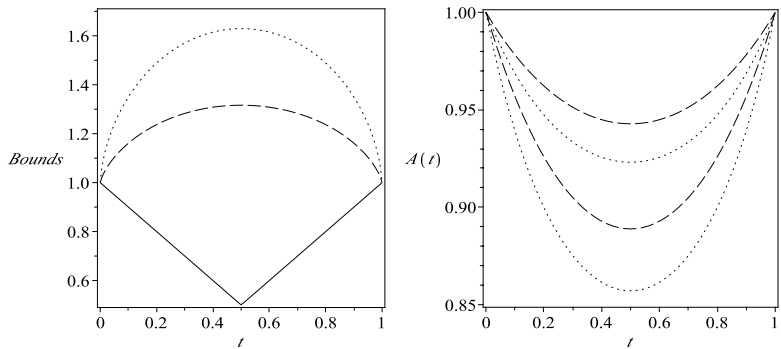}

\caption{Left panel: graph of the bounds $A_W^\mathrm{P}$ (top curve),
$A_W^\mathrm{CFG}$ (middle curve) and $A_M^\mathrm{P} = A_M^\mathrm{CFG}$
(bottom curve). Right panel: graph of $A_C^\mathrm{P}$ (dashed) and
$A_C^\mathrm{CFG}$ (dotted) for the Farlie--Gumbel--Morgenstern copula
with $\theta= 1/2$ (upper curves) and $1$ (lower curves).}
\label{fig:1}
\end{figure}

The bounds $A_W^\mathrm{P}$, $A_W^\mathrm{CFG}$ and $A_M^\mathrm{P} =
A_M^\mathrm{CFG}$ are depicted in the left panel of Figure \ref{fig:1}.
As a
further example, consider the Farlie--Gumbel--Morgenstern copula with
parameter $\theta\in[-1,1]$, defined for all $u, v \in[0,1]$ by
$C_\theta(u,v) = uv + \theta uv(1-u)(1-v)$. Condition \eqref{eq:7}
is met if $\theta\ge0$ and it is easy to check that for all $t \in[0,1]$,
%
%
\begin{equation}\label{eq:8}
A^\mathrm{P}_{\theta} (t) = \frac{2t^2-2t-4}{2t^2-2t-4 + (3t^2
-3t)\theta}
,\qquad
A^\mathrm{CFG}_{\theta} (t) = \biggl( \frac{2}{2+t-t^2} \biggr)^\theta.
\end{equation}
These functions are graphed in the right panel of Figure \ref{fig:1}.

Invoking Proposition \ref{prop:2}, we can proceed as in Appendix \hyperref[appa]{A} to
show the convergence of the goodness-of-fit process in the case of LTD
copulas, whence the following result. The parametric bootstrap
algorithm described in Section \ref{sec:3} also applies \textit{mutatis
mutandis} and remains valid under such $H_0^*$.

\begin{prop}
\label{prop:4}
Assume $H_0^*$ holds, that is, $C$ is an LTD copula such that $C =
C_{\theta_0}$ for some $\theta_0 \in\mathcal{O}$. Let $\mathcal
{A}^\mathrm{P} = \{ A_C^\mathrm{P}\dvtx  C \in\mathcal{C} \}$ and $\mathcal
{A}^\mathrm{CFG} = \{ A_C^\mathrm{CFG}\dvtx  C \in\mathcal{C} \}$.

\begin{enumerate}[(a)]
\item[(a)] If $\mathcal{A}^\mathrm{P}$ meets conditions \textup
{(A)}--\textup{(C)} and $(\mathbb{A}_n^{\mathrm{P}}, \Theta_n)$
converges to a Gaussian limit
$(\mathbb{A}^{\mathrm{P}}_C ,\Theta)$, then $\mathbb{A}_{n,\theta
_n} \leadsto\mathbb{A}^{\mathrm{P}}_C - {\dot
A}_{\theta_0}^\top\Theta$ as $n \to\infty$ in $\mathcal{C} [0,1]$.
\item[(b)] If $\mathcal{A}^\mathrm{CFG}$ meets conditions
\textup{(A)}--\textup{(C)} and $(\mathbb{A}_n^{\mathrm{CFG}}, \Theta_n)$
converges to a Gaussian limit
$(\mathbb{A}^{\mathrm{CFG}}_C ,\Theta)$, then $\mathbb{A}_{n,\theta
_n} \leadsto\mathbb{A}^{\mathrm{CFG}}_C - {\dot
A}_{\theta_0}^\top\Theta$ as $n \to\infty$ in $\mathcal{C} [0,1]$.
\end{enumerate}
\end{prop}

\section{Consistency of the test}
\label{sec:5}

Suppose that $C \notin\mathcal{C}$ is an LTD copula and that the
hypothesis $H_0^*\dvtx  C \in\mathcal{C}$ is being tested with the Cram\'
er--von Mises statistic $S_n$. Let $A_n$ denote either $A_n^{\mathrm
{P}}$ or $A_n^{\mathrm{CFG}}$
and let $A$ stand for $A_C^\mathrm{P}$ or $A_C^\mathrm{CFG}$, as the
case may
be. Further, assume that $\theta_n$ is a consistent, rank-based
estimator of some $\theta^* \in\mathcal{O}$. The test based on $S_n$
is then consistent, provided that $A \neq A_{\theta^*}$.

To see this, decompose the process $\mathbb{A}_{n,\theta_n}$ as
%
%
\begin{equation}
\label{eq:9}
\sqrt{n} ( A_n - A_{\theta_n} ) = \sqrt{n} (A_n -A) - \sqrt
{n}
(A_{\theta_n} - A_{\theta^*} ) + \sqrt{n} (A - A_{\theta^*}).
\end{equation}
Assume conditions (A)--(C) hold for $\mathcal{A} = \mathcal{A}^\mathrm{P}$
or $\mathcal{A}^\mathrm{CFG}$ and that as $n \to\infty$, $(\sqrt{n}
(A_n - A), \sqrt{n} (\theta_n - \theta^*)) \rightsquigarrow
(\mathbb{A},
\Theta^*)$ to a Gaussian limit, where $\mathbb{A}$ stands for either
$\mathbb{A}^{\mathrm{P}}$ or
$\mathbb{A}^{\mathrm{CFG}}$. We can then proceed exactly as in
Appendix \hyperref[appa]{A} to see that as
$n \to\infty$, $\sqrt{n} (A_n - A) - \sqrt{n} (A_{\theta_n} -
A_{\theta^*} ) \rightsquigarrow\mathbb{A}- {\dot A}_{\theta^*}
^\top\Theta
^*$. If $A \neq A_{\theta^*}$, then $\sup_{t \in[0,1]} \sqrt{n}
|A(t) - A_{\theta^*}(t)| \to\infty$ and hence, for every $\epsilon>0$,
\[
\lim_{n\to\infty}\Pr(S_n > \epsilon) = 1.
\]

In particular, the test based on $S_n$ is consistent whenever $C$ is an
extreme-value copula and the hypothesized family $\mathcal{C}$ also
consists of extreme-value copulas. However, consistency may fail
otherwise, for it may happen that $A = A_{\theta^*}$, even if
$H_0^\ast
$ is false.

To illustrate this point, consider the functions $A^\mathrm{P}_{\theta}$
and $A^\mathrm{CFG}_{\theta}$ given in \eqref{eq:8}. As the latter are
convex, they can be used to generate new families of extreme-value
copulas, which may be called the FGM--P and FGM--CFG families.

Now, suppose that $C$ is the Farlie--Gumbel--Morgenstern copula with
parameter $\theta> 0$ and that the statistic $S_n$ is used to test
$H_0\dvtx  A \in\mathcal{A}$ when:
\begin{longlist}[(a)]
\item[(a)] $\mathcal{A}$ is the Gumbel--Hougaard family of copulas;
\item[(b)] $\mathcal{A}$ is the FGM--CFG family of extreme-value copulas.
\end{longlist}
In case (a), the tests based on $A_n^{\mathrm{P}}$ and $A_n^{\mathrm
{CFG}}$ would be consistent
because $A^\mathrm{P}_\theta$ and $A^\mathrm{CFG}_\theta$ both differ from
the Pickands dependence function of the Gumbel--Hougaard given in
\eqref{eq:3}. In case (b), the test based on $A_n^{\mathrm{P}}$
would also be
consistent because $A^\mathrm{P}_\theta\neq A^\mathrm{CFG}_{\theta^*}$. The
test based on $A_n^{\mathrm{CFG}}$ may fail to be consistent, however,
given that
$A^\mathrm{CFG}_\theta$ coincides with the Pickands dependence function of
the FGM--CFG family. Consistency of the test would then depend on the
behavior of $\theta_n$.

Suppose, for instance, that $\theta$ is estimated by inversion of
Kendall's tau. As $n \to\infty$, $\theta_n$ would approach $2
\theta
/9$, which is the population value of this dependence measure for the
FGM copula. For the FGM--CFG family, however, Kendall's tau is $7\theta
/10 + \theta^2/30$, which coincides with $2 \theta/9$ only when
$\theta= 0$, that is, at independence where the difference between the
two models is immaterial. Therefore, the test based on $A_n^{\mathrm
{CFG}}$ would be
consistent in this case, provided that $\theta$ is estimated by
inversion of Kendall's tau. A similar conclusion would be reached for
inversion of Spearman's rho and maximum pseudo-likelihood estimation.

\section{Power study}
\label{sec:6}

Equation \eqref{eq:9} and the accompanying discussion suggest that
just as for consistency, the power of the test based on $S_n$ depends
on how different $A = A_C^\mathrm{P}$ or $A_C^\mathrm{CFG}$ is from its
parametric estimate $A_{\theta^*}$ under $H_0$. This issue is
investigated graphically in Section \ref{sec:6.1} and via simulations
in Sections~\ref{sec:6.2} and \ref{sec:6.3}.

\subsection{General considerations}
\label{sec:6.1}

Consider the following three sets of LTD copula families.
\begin{description}
\item[Group I:] \textit{Symmetric extreme-value copulas}: the
Gumbel--Hougaard (GH), Galambos (GA), H\"usler--Reiss (HR) and Student
extreme-value (t-EV) copula with four degrees of freedom.
\item[Group II:] \textit{Symmetric non-extreme-value copulas}: the
Clayton (C), Frank (F), Normal (N) and Plackett (P).
\item[Group III:] \textit{Asymmetric extreme-value copulas}: asymmetric
versions of the Gumbel--Hougaard (a-GH), Galambos (a-GA), H\"
usler--Reiss (a-HR), and Student extreme-value (a-t-EV) copula with
four degrees of freedom.
\end{description}

\begin{figure}

\includegraphics{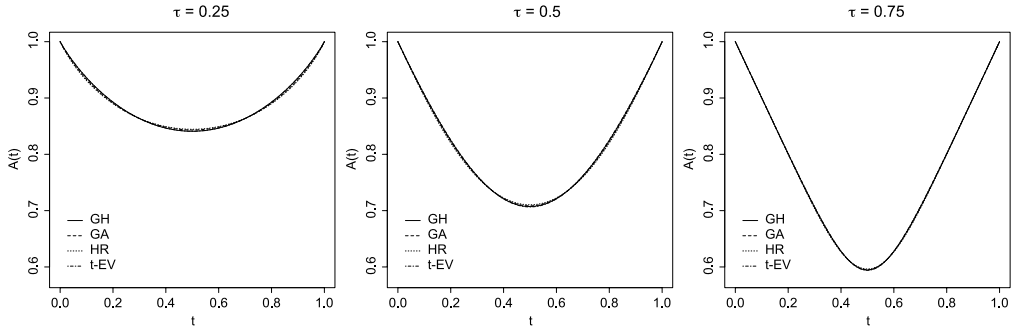}

\caption{
Pickands dependence functions of the Gumbel--Hougaard, Galambos, H\"
usler--Reiss and t-EV copulas when $\tau= 0.25$, $\tau= 0.50$ and
$\tau= 0.75$.}\label{fig:2}\vspace*{-2pt}
\end{figure}

Figure \ref{fig:2} shows the Pickands dependence functions of the
copulas in Group I when $\tau= 0.25$, $0.50$, $0.75$. Although the
curves are not identical, they are very similar. When the statistic
$S_n$ is used to distinguish between these models, therefore, the test
will be consistent, but can be expected to have little power, even in
moderate sample sizes.

%

\begin{figure}[b]\vspace*{-2pt}

\includegraphics{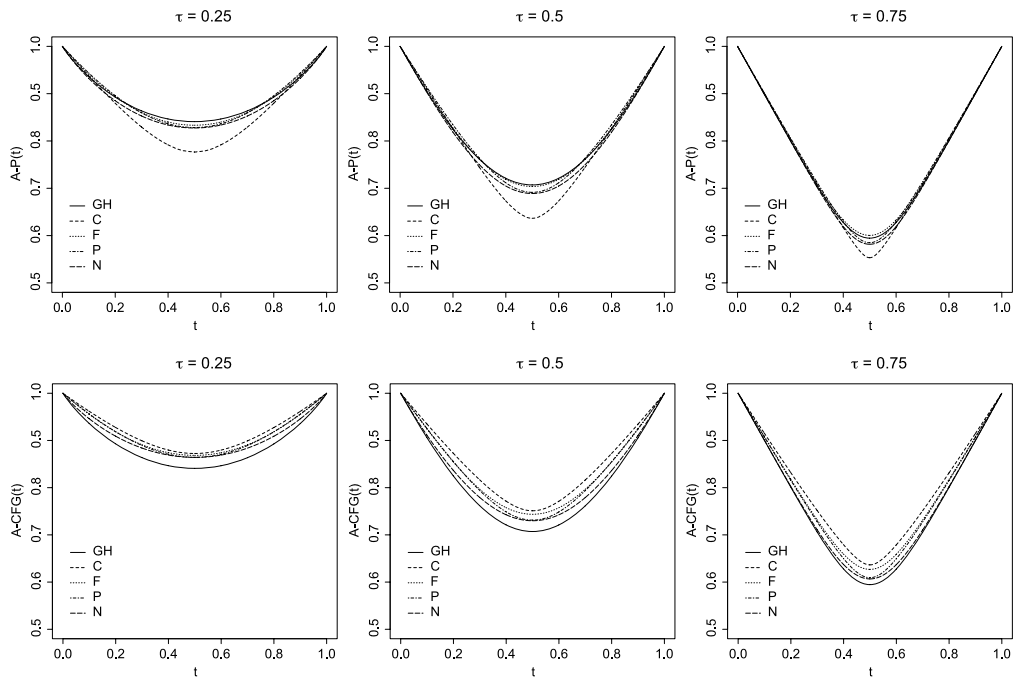}

\caption{Plots of $A_C^\mathrm{P}$ (top) and $A_C^\mathrm{CFG}$ (bottom) when
$C$ is
the Gumbel--Hougaard (GH), Clayton (C), Frank (F), Normal (N) and
Plackett (P) copula with $\tau= 0.25$ (left), $\tau= 0.50$ (middle)
and $\tau= 0.75$ (right).}
\label{fig:3}
\end{figure}

In Figure \ref{fig:3}, the functions $A_C^\mathrm{P}$ and $A_C^\mathrm{CFG}$
are plotted for the copulas in Group II and the same values of tau. For
comparison purposes, the curve corresponding to the Gumbel--Hougaard
copula is added. Here, the differences between the curves are much more
pronounced. Thus, the power of the test based on $S_n$ may be expected
to rise quickly (and be approximately the same) if the copula family
under $H_0$ is from Group I.

\begin{figure}

\includegraphics{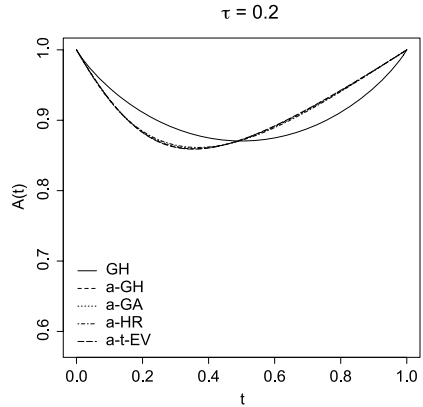}

\caption{Pickands dependence functions for the Gumbel--Hougaard copula and four
asymmetric extreme-value copulas with $\tau= 0.20$: the asymmetric
Gumbel--Hougaard (a-GH), Galambos (a-GA), H\"usler--Reiss (a-HR) and
t-EV (a-t-EV) with four degrees of freedom.}
\label{fig:4}
\end{figure}

Figure \ref{fig:4} shows the Pickands dependence functions of the
copulas in Group III. These copulas were derived using Khoudraji's
device \cite{Kho95,GenGhoRiv98,Lieb08}, which transforms any
symmetric copula $C_\theta$ into a non-exchangeable model via the formula
\[
C_{\lambda, \kappa, \theta}(u,v) = u^{1-\lambda} v^{1-\kappa
}C_\theta
(u^\lambda,v^\kappa)
\]
for all $u, v \in[0,1]$ and arbitrary choices of $\lambda\neq\kappa
\in(0,1)$. Furthermore, if $C_\theta$ is an extreme-value copula with
Pickands dependence function $A_\theta$, then $C_{\lambda, \kappa,
\theta}$ is also an extreme-value copula. Its Pickands dependence
function is given, at all $t \in[0,1]$, by
\[
A_{\lambda, \kappa, \theta} (t) = (1-\kappa)t + (1-\lambda)(1-t) +
\{
\kappa t + \lambda(1-t) \} A_\theta\biggl\{ \frac{\kappa t}{\kappa t +
\lambda(1-t)} \biggr\}.
\]

%

Note that the dependence in $C_{\lambda, \kappa, \theta}$ is limited
since, by the Fr\'echet--Hoeffding inequality,
\[
C_{\lambda, \kappa, \theta} (u,v) \le u^{1-\lambda} v^{1-\kappa}
\min
(u^\lambda, v^\kappa) = \min(uv^{1-\kappa}, vu^{1-\lambda}).
\]
As the right-hand term is the Marshall--Olkin copula MO$_{\lambda
,\kappa}$, Example 5.5 in \cite{Nelsen2006}, implies that
\[
\tau(C_{\lambda, \kappa, \theta}) \le\tau(\mathrm{MO}_{\lambda
,\kappa})
= \frac{\kappa\lambda}{\kappa+ \lambda- \kappa\lambda} .
\]

In the present study, the values $\lambda= 0.3$, $\kappa=0.8$ were
used and, hence, $\tau(C_{\lambda, \kappa, \theta})$ could not exceed
$0.279$. For each choice of copula family $C_\theta$ in Group III, the
parameter $\theta$ was set to make Kendall's tau equal to $0.20$.

Figure \ref{fig:4} shows that the Pickands dependence functions of the
copulas in Group III are very similar, though distinct. They are,
however, easily distinguished from their symmetric counterparts with
the same value of tau. Thus, although these extreme-value copulas would
be difficult to tell apart on the basis of $S_n$ in moderate samples,
the test may still be reasonably powerful against copulas in Group I.

\subsection{Monte Carlo study}
\label{sec:6.2}

The observations in Section \ref{sec:6.1} were confirmed through
simulations. To this end, $1000$ random samples of size $n = 300$ were
generated from 28 different copulas, $C$, corresponding to the
following scenarios:
\begin{longlist}[(a)]
\item[(a)] $C$ belongs to Group I or II and $\tau(C) \in\{ 0.25,
0.50, 0.75\}$;
\item[(b)] $C$ belongs to Group III and $\tau(C) = 0.20$.
\end{longlist}
The statistics $S_n^{\mathrm{P}}$ and $S_n^{\mathrm{CFG}}$ were
computed for each data set. Four
hypotheses of the form $H_0\dvtx  A \in\mathcal{A}$ were then tested. The
choices for $\mathcal{A}$ were the families of Pickands dependence
functions for extreme-value copulas in Group I.

All tests were carried out at the 5\% level. Each $P$-value was
computed on the basis of $N=1000$ parametric bootstrap samples. For
comparison purposes, goodness of fit was also checked with the general
purpose statistic
\[
T_n = \sum_{i=1}^n | C_n (U_i, V_i) - C_{\theta_n}(U_i, V_i) |^2.
\]
This particular test statistic was chosen because of its good overall
performance in the large scale simulation studies of Berg \cite{Ber09} and
Genest \textit{et al.} \cite{GenRemBea09}.

Tables \ref{tab:1}--\ref{tab:4} report the percentages of rejection of
the four null hypotheses under each scenario. Although this made little
difference, these results are for the end-point-corrected versions of
$A_n^{\mathrm{P}}$ and $A_n^{\mathrm{CFG}}$, defined for all $t \in
[0,1]$ by
\[
1/A_{n,c}^{\mathrm{P}}(t) = 1/A_n^{\mathrm{P}}(t) -(1-t)\{
1/A_n^{\mathrm{P}}(0)-1\} - t \{ 1/A_n^{\mathrm{P}}(1)-1\}
\]
and
\[
\log A_{n,c}^{\mathrm{CFG}}(t) = \log\{ A_n^{\mathrm{CFG}}(t)\}
-(1-t) \log\{ A_n^{\mathrm{CFG}}(0)\} - t
\log\{ A_n^{\mathrm{CFG}}(1)\}.
\]

Before commenting on the results, note that for copulas in Groups I and
II, the real-valued dependence parameter of each data set was estimated
by inversion of Kendall's tau; its implementation relied on the
numerical approximation technique of Kojadinovic and Yan \cite{KojYan09}. For copulas in
Group III, which involve several parameters, maximum pseudo-likelihood
estimation was used \cite{GenGhoRiv95,ShiLou95}.

%
\begin{table}
\caption{Percentage of rejection of $H_0$ for copulas in Group I when $n=300$}
\label{tab:1}
\begin{tabular*}{\textwidth}{@{\extracolsep{\fill}}llccccccccccc@{}}
\hline
$H_0$ & True & \multicolumn{3}{c}{$\tau= 0.25$} & &
\multicolumn{3}{c}{$\tau= 0.50$} & & \multicolumn{3}{c@{}}{$\tau=
0.75$} \\[-6pt]
& & \multicolumn{3}{c}{\hrulefill} & &
\multicolumn{3}{c}{\hrulefill} & & \multicolumn{3}{c@{}}{\hrulefill} \\
& & $T_n$ & $S_n^{\mathrm{P}}$ & $S_n^{\mathrm{CFG}}$ & & $T_n$ &
\multicolumn{1}{c}{$S_n^{\mathrm{P}}$\phantom{0}} & \multicolumn{1}{c}{$S_n^{\mathrm{CFG}}$} & & $T_n$ &
\multicolumn{1}{c}{$S_n^{\mathrm{P}}$\phantom{0}} & \multicolumn{1}{c@{}}{$S_n^{\mathrm{CFG}}$} \\
\hline
GH & GH & \textbf{4.2} & \textbf{3.8} & \textbf{4.0} & & \textbf{4.0} &
\textbf{4.8} & \textbf{3.6} & & \textbf{4.2} & \textbf{5.3} & \textbf{5.5} \\
& GA & 4.8 & 4.3 & 4.2 & & 4.4 & 3.8 & 3.8 & & 4.8 & 4.3 & 4.3 \\
& HR & 4.8 & 4.2 & 4.0 & & 5.4 & 3.4 & 3.9 & & 3.7 & 3.1 & 1.7 \\
& t-EV & 4.2 & 3.8 & 4.5 & & 5.1 & 5.5 & 6.4 & & 4.8 & 7.5 & 8.9
\\[6pt]
GA
& GH & 4.5 & 4.7 & 3.9 & & 4.0 & 5.8 & 4.7 & & 4.4 & 5.6 & 6.8 \\
& GA & \textbf{4.3} & \textbf{4.6} & \textbf{4.0} & & \textbf{5.5} &
\textbf{3.9} &
\textbf{4.8} & & \textbf{4.3} & \textbf{4.7} & \textbf{4.6} \\
& HR & 4.6 & 4.8 & 4.2 & & 5.0 & 3.4 & 3.4 & & 3.7 & 3.7 & 1.8 \\
& t-EV & 4.6 & 4.7 & 4.4 & & 5.3 & 8.0 & 7.1 & & 5.7 & 8.2 &
10.9\phantom{0} \\[6pt]
HR
& GH & 4.6 & 6.4 & 4.4 & & 4.3 & 9.6 & 7.5 & & 4.5 & 9.6 & 15.7\phantom{0}
\\
& GA & 4.3 & 5.4 & 4.5 & & 5.1 & 6.6 & 7.2 & & 5.1 & 8.4 & 11.7\phantom{0}
\\
& HR & \textbf{4.9} & \textbf{5.2} & \textbf{4.2} & & \textbf{5.3} &
\textbf{4.3} &
\textbf{3.9} & & \textbf{4.0} & \textbf{4.3} & \textbf{3.3} \\
& t-EV & 4.6 & 5.9 & 4.8 & & 5.8 & 13.7\phantom{0} & 11.5\phantom{0} & & 6.6 & 14.9\phantom{0} &
29.3\phantom{0} \\[6pt]
t-EV
& GH & 4.2 & 3.4 & 4.0 & & 4.1 & 3.9 & 2.9 & & 4.0 & 3.3 & 2.4 \\
& GA & 4.1 & 4.3 & 4.4 & & 4.8 & 3.4 & 3.9 & & 4.6 & 3.0 & 1.7 \\
& HR & 4.7 & 4.1 & 4.4 & & 5.4 & 3.2 & 3.4 & & 3.8 & 2.2 & 1.3 \\
& t-EV & \textbf{4.6} & \textbf{3.7} & \textbf{4.2} & & \textbf{4.7} &
\textbf{4.8}
& \textbf{5.2} & & \textbf{4.1} & \textbf{4.3} & \textbf{4.7} \\
\hline
\end{tabular*}\vspace*{-3pt}
\end{table}
%
\begin{table}
\caption{Percentage of rejection of $H_0$ for copulas in Group I when $n=1000$}
\label{tab:2}
\begin{tabular*}{\textwidth}{@{\extracolsep{\fill}}llccccccccccc@{}}
\hline
$H_0$ & True & \multicolumn{3}{c}{$\tau= 0.25$} & &
\multicolumn{3}{c}{$\tau= 0.50$} & & \multicolumn{3}{c@{}}{$\tau=
0.75$} \\[-6pt]
& & \multicolumn{3}{c}{\hrulefill} & &
\multicolumn{3}{c}{\hrulefill} & & \multicolumn{3}{c@{}}{\hrulefill} \\
& & $T_n$ & \multicolumn{1}{c}{$S_n^{\mathrm{P}}$\phantom{0}} & \multicolumn{1}{c}{$S_n^{\mathrm{CFG}}$} & & $T_n$ &
\multicolumn{1}{c}{$S_n^{\mathrm{P}}$\phantom{0}} & \multicolumn{1}{c}{$S_n^{\mathrm{CFG}}$} & & $T_n$ &
\multicolumn{1}{c}{$S_n^{\mathrm{P}}$\phantom{0}} & \multicolumn{1}{c@{}}{$S_n^{\mathrm{CFG}}$} \\
\hline
GH
& GH & \textbf{4.9} & \textbf{4.7} & \textbf{5.3} & & \textbf{4.0} &
\textbf{5.9} &
\textbf{6.0} & & \textbf{3.8} & \textbf{5.2} & \textbf{5.4} \\
& GA & 5.1 & 5.8 & 4.0 & & 5.9 & 4.4 & 5.1 & & 4.8 & 3.8 & 4.2 \\
& HR & 5.1 & 6.3 & 6.3 & & 5.1 & 6.3 & 9.0 & & 3.4 & 3.5 & 9.2 \\
& t-EV & 5.4 & 4.4 & 5.4 & & 6.1 & 6.2 & 6.9 & & 5.4 & 9.8 &
15.9\phantom{0} \\[6pt]
GA
& GH & 5.2 & 7.4 & 6.1 & & 4.4 & 8.1 & 8.4 & & 4.3 & 5.6 & 7.3 \\
& GA & \textbf{5.0} & \textbf{5.6} & \textbf{4.0} & & \textbf{5.4} &
\textbf{5.1} &
\textbf{5.4} & & \textbf{4.8} & \textbf{4.4} & \textbf{5.2} \\
& HR & 4.4 & 5.0 & 5.2 & & 4.5 & 4.5 & 6.2 & & 3.5 & 3.1 & 6.2 \\
& t-EV & 6.1 & 6.9 & 6.6 & & 6.7 & 9.4 & 12.7\phantom{0} & & 5.5 & 12.8\phantom{0} &
23.1\phantom{0} \\[6pt]
HR
& GH & 6.2 & 10.6\phantom{0} & 8.6 & & 5.1 & 17.6\phantom{0} & 17.8\phantom{0} & & 5.5 & 18.1\phantom{0} &
40.2\phantom{0} \\
& GA & 5.4 & 6.6 & 4.1 & & 5.6 & 8.1 & 8.8 & & 5.5 & 12.7\phantom{0} & 23.4\phantom{0}
\\
& HR & \textbf{4.6} & \textbf{5.9} & \textbf{5.5} & & \textbf{4.2} &
\textbf{4.9} &
\textbf{5.1} & & \textbf{3.4} & \textbf{4.7} & \textbf{5.6} \\
& t-EV & 6.6 & 10.1\phantom{0} & 8.2 & & 8.2 & 27.0\phantom{0} & 34.4\phantom{0} & & 6.5 & 45.2\phantom{0} &
81.7\phantom{0} \\[6pt]
t-EV
& GH & 4.7 & 4.7 & 5.3 & & 4.4 & 4.4 & 5.7 & & 4.0 & 3.4 & 3.0 \\
& GA & 4.8 & 5.6 & 4.2 & & 5.6 & 4.8 & 6.0 & & 4.8 & 3.0 & 3.7 \\
& HR & 5.3 & 6.4 & 6.1 & & 5.5 & 8.5 & 12.3\phantom{0} & & 4.3 & 3.9 & 28.7\phantom{0}
\\
& t-EV & \textbf{5.1} & \textbf{4.5} & \textbf{5.5} & & \textbf{5.6} &
\textbf{4.8}
& \textbf{5.3} & & \textbf{5.2} & \textbf{4.5} & \textbf{4.7} \\
\hline
\end{tabular*}

\end{table}
%

\subsection{Results}
\label{sec:6.3}

It is clear from Table \ref{tab:1} that when $n=300$, the tests based
on $T_n$, $S_n^{\mathrm{P}}$ and $S_n^{\mathrm{CFG}}$ cannot
distinguish between copulas in
Group I. When $\tau= 0.25$, all rejection rates are within sampling
error from the nominal level. There are only small signs of improvement
as $\tau$ rises to $0.50$ and $0.75$. The best scores are obtained when
testing for the H\"usler--Reiss model with $S_n^{\mathrm{CFG}}$ when
$\tau= 0.75$.
Globally, there is little to choose between the tests.

%
\begin{table}
\caption{Percentage of rejection of $H_0$ for copulas in Group II when $n=300$}
\label{tab:3}
\begin{tabular*}{\textwidth}{@{\extracolsep{\fill}}lld{2.1}d{2.1}d{2.1}d{3.1}d{3.1}d{3.1}d{3.1}d{3.1}d{3.1}@{}}
\hline
$H_0$ & True & \multicolumn{3}{c}{$\tau= 0.25$}  &
\multicolumn{3}{c}{$\tau= 0.50$}  & \multicolumn{3}{c@{}}{$\tau=
0.75$} \\[-6pt]
& & \multicolumn{3}{c}{\hrulefill}  &
\multicolumn{3}{c}{\hrulefill} & \multicolumn{3}{c@{}}{\hrulefill} \\
& & \multicolumn{1}{c}{$T_n$} & \multicolumn{1}{c}{$S_n^{\mathrm{P}}$} & \multicolumn{1}{c}{$S_n^{\mathrm{CFG}}$} &
\multicolumn{1}{c}{$T_n$} &
\multicolumn{1}{c}{$S_n^{\mathrm{P}}$} & \multicolumn{1}{c}{$S_n^{\mathrm{CFG}}$} & \multicolumn{1}{c}{$T_n$} &
\multicolumn{1}{c}{$S_n^{\mathrm{P}}$} & \multicolumn{1}{c@{}}{$S_n^{\mathrm{CFG}}$} \\
\hline
GH & C & 98.8 & 99.5 & 82.1 & 100.0 & 100.0 & 100.0 &100.0 & 100.0 &
100.0 \\
& F & 36.6 & 11.0 & 48.0 & 82.0 & 7.1 & 100.0 & 92.1 & 27.2 & 100.0
\\
& N & 26.9 & 21.9 & 21.8 & 43.5 & 44.9 & 66.9 & 37.5 & 18.7 & 82.3
\\
& P & 34.3 & 17.3 & 43.4 & 68.0 & 44.6 & 98.6 & 65.0 & 71.6 & 100.0
\\[6pt]
GA & C & 98.9 & 99.7 & 84.0 & 100.0 & 100.0 & 100.0 &100.0 & 100.0 &
100.0 \\
& F & 39.8 & 15.1 & 50.1 & 83.4 & 10.2 & 100.0 & 92.1 & 29.9 &
100.0 \\
& N & 28.1 & 25.7 & 21.9 & 44.0 & 49.0 & 69.5 & 37.4 & 21.5 & 83.2
\\
& P & 37.7 & 23.4 & 45.0 & 70.8 & 57.1 & 99.0 & 65.7 & 76.7 & 100.0
\\[6pt]
HR & C & 99.1 & 99.9 & 84.5 & 100.0 & 100.0 & 100.0 &100.0 & 100.0 &
100.0 \\
& F & 42.3 & 18.8 & 52.5 & 85.2 & 18.9 & 100.0 & 93.7 & 42.1 &
100.0 \\
& N & 28.3 & 29.0 & 22.5 & 46.0 & 55.7 & 73.3 & 38.9 & 34.8 & 89.8
\\
& P & 41.1 & 28.8 & 48.3 & 75.1 & 74.5 & 99.5 & 73.1 & 92.1 & 100.0
\\[6pt]
t-EV & C & 98.6 & 99.5 & 82.6 & 100.0 & 100.0 & 100.0 &100.0 & 100.0
& 100.0 \\
& F & 36.7 & 11.1 & 48.3 & 81.3 & 5.0 & 100.0 & 90.9 & 18.9 & 100.0
\\
& N & 26.5 & 21.8 & 21.7 & 43.5 & 42.7 & 66.2 & 36.9 & 10.6 & 74.2
\\
& P & 34.8 & 17.2 & 43.7 & 67.7 & 35.7 & 98.1 & 62.1 & 53.2 & 99.8
\\
\hline
\end{tabular*}
\end{table}
%
%
\begin{table}
\caption{Percentage of rejection of $H_0$ for copulas in Group III
when $n=300$}
\label{tab:4}
\begin{tabular*}{250pt}{@{\extracolsep{\fill}}lld{2.1}d{2.1}d{2.1}@{}}
\hline
\multicolumn{1}{@{}l}{$H_0$} & True & \multicolumn{1}{c}{$T_n$} & \multicolumn{1}{c}{$S_n^{\mathrm{P}}$} & \multicolumn{1}{c@{}}{$S_n^{\mathrm{CFG}}$} \\
\hline
GH & a-GH & 32.7 & 40.9 & 86.5 \\
& a-GA & 33.5 & 42.8 & 86.7 \\
& a-HR & 28.4 & 37.5 & 83.5 \\
& a-t-EV & 33.1 & 41.4 & 88.6 \\[6pt]
 GA& a-GH & 33.4 & 40.8 & 89.2 \\
& a-GA & 34.0 & 42.3 & 89.3 \\
& a-HR & 28.4 & 38.1 & 86.5 \\
& a-t-EV & 32.7 & 40.6 & 90.5 \\[6pt]
HR & a-GH & 36.2 & 37.5 & 93.3 \\
& a-GA & 31.7 & 39.1 & 89.7 \\
& a-HR & 32.6 & 40.9 & 90.3 \\
& a-t-EV & 40.5 & 42.8 & 92.3 \\[6pt]
t-EV
& a-GH & 32.0 & 41.2 & 87.1 \\
& a-GA & 33.2 & 43.3 & 87.8 \\
& a-HR & 27.3 & 38.4 & 83.6 \\
& a-t-EV & 31.3 & 40.7 & 88.7 \\
\hline
\end{tabular*}
\end{table}

Table \ref{tab:2} shows what happens when $n=1000$. Power is on the
rise, especially when $\tau= 0.75$. In the latter case, it seems
preferable to base the test on $S_n^{\mathrm{CFG}}$ rather than on
$S_n^{\mathrm{P}}$ -- both do
better than the test based on $T_n$. Overall, the results remain
disappointingly low, except when testing for the H\"usler--Reiss model
with $\tau\ge0.50$.

These observations are in line with Figure \ref{fig:2}, which shows
striking similarities between the Gumbel--Hougaard, Galambos, H\"
usler--Reiss and t-EV copula with four degrees of freedom. While
$S_n^{\mathrm{P}}$
and $S_n^{\mathrm{CFG}}$ still have difficulty telling them apart when
the sample
size is $1000$, their power eventually rises when $n \to\infty$, as
explained in Section \ref{sec:5}. To illustrate this point, samples of
various sizes were generated from the Gumbel--Hougaard copula with
$\tau= 0.50$ and the statistic $S_n^{\mathrm{CFG}}$ was used to test
for the
Galambos family. The following results, based on $1000$ repetitions and
$N=1000$ bootstrap samples, give an idea of the sample sizes needed to
differentiate models in Group I:

\begin{center}
\begin{tabular}{@{}lllll@{}}
\hline
Sample size $n$ & 5\,000 & 10\,000 & 20\,000 &
40\,000 \\
\hline
Percentage of rejection of $H_0$ & 10.8 & 22.6 & 60.2 & 97.3\\
\hline
\end{tabular}\vspace*{6pt}
\end{center}

Returning to the case $n = 300$, we can see from Table \ref{tab:3} that
the test based on $S_n^{\mathrm{CFG}}$ is quite good at detecting
non-extreme-value
LTD alternatives from Group II. Its power is higher than those of
$S_n^{\mathrm{P}}$
and $T_n$, except when the data are generated from the Clayton or the
Normal copula with $\tau= 0.25$. Interestingly, the general purpose
test based on $T_n$ is often second best. The statistic $S_n^{\mathrm
{P}}$ has the
edge only for the Clayton when $\tau= 0.25$; it does very poorly
against the Frank, and against the Normal when $\tau= 0.75$.

These results are in close agreement with the plots displayed in
Figure \ref{fig:3}. Consider, for instance, the case where
$S_n^{\mathrm{P}}$ is
used to test for the Gumbel--Hougaard copula from weakly dependent data
($\tau= 0.25)$. From Table \ref{tab:3}, the alternatives can be ranked
as follows in decreasing order of power:
\[
\mbox{Clayton  }\succ\mbox{ Normal  }\succ\mbox{ Plackett }\succ\mbox{ Frank}.
\]
Looking at Figure \ref{fig:3}, we find that this ordering is concordant
with the overall degree of dissimilarity between $A_C^\mathrm{P}$ and $A$.
In this case, as in others, it is found that at fixed sample size,
curves that look alike are harder to distinguish than others.

Finally, Table \ref{tab:4} shows that the statistic $S_n^{\mathrm
{CFG}}$ is much
better than the other two at detecting asymmetric extreme-value
alternatives. The overall good performance of this test is consistent
with evidence from \cite{GenSeg09} that $A_n^{\mathrm{CFG}}$ is
generally a better
nonparametric estimator of the Pickands dependence function than
$A_n^{\mathrm{P}}$.
When the margins are known, this phenomenon is well documented; see,
for example, \cite{CapFouGen97,HallTajvidi2000} or \cite
{Jimenez2001}.

\section{Conclusion}
\label{sec:7}

Copula models are now common. As illustrated, for instance, by Ben Ghorbal \textit{et al.} \cite
{BenGenNes09}, so are situations in which the dependence structure of a
random pair $(X,Y)$ is well represented by an extreme-value copula,
even though $X$ and $Y$ themselves do not necessarily exhibit
extreme-value behavior. In such cases, the statistics considered here
can be used to test the goodness of fit of specific parametric copula
families of the form \eqref{eq:2} such as the Gumbel--Hougaard,
Galambos, H\"usler--Reiss or Student extreme-value copula.

Theoretical and empirical evidence presented here shows that the
nonparametric tests based on the Cram\'er--von Mises statistic $S_n$
are generally consistent and that they are an effective tool for
distinguishing between symmetric and asymmetric extreme-value copulas,
as well as for detecting other left-tail decreasing (LTD) dependence structures.

Except in the presence of massive data, however, it seems very
difficult to discriminate between extreme-value copulas whose Pickands
dependence functions are close. This may come as something of a
disappointment, but, on reflection, we may wonder whether, in the light
of Figure \ref{fig:2}, there is any \textit{practical} difference
between, say, the Gumbel--Hougaard and the Galambos copula when they
have the same value of Kendall's tau.

For example, many studies have concluded that a Gumbel--Hougaard copula
structure is adequate for the insurance data mentioned in the
Introduction; see, for instance, \cite{FreVal98,GenGhoRiv98,CheFan05,Denuit2006,DupuisJones2006,GenQueRem06b} or
\cite{KojYan09c}. In these papers, comparisons were
made between the Gumbel--Hougaard model and non-extreme-value copulas
that were either Archimedean or meta-elliptical.

As Ben Ghorbal \textit{et al.} \cite{BenGenNes09} conclude that the data exhibit extreme-value
dependence, it may be worth comparing the Gumbel--Hougaard structure
with other extreme-value copulas from Groups I and III. This is done in
Table \ref{tab:5} using the statistics $S_n^{\mathrm{P}}$ and
$S_n^{\mathrm{CFG}}$ and the
inversion of Kendall's tau to estimate $\theta$. Because the test is
yet to be adapted to the case of censoring, the analysis ignored the 34
claims for which the policy limit was reached. Each $P$-value in the
table is based on $N = 2500$ bootstrap samples. Given the comparatively
small sample size, $n = 1466$, it is little wonder that no
model is rejected at the 5\% level.

\begin{figure}[b]

\includegraphics{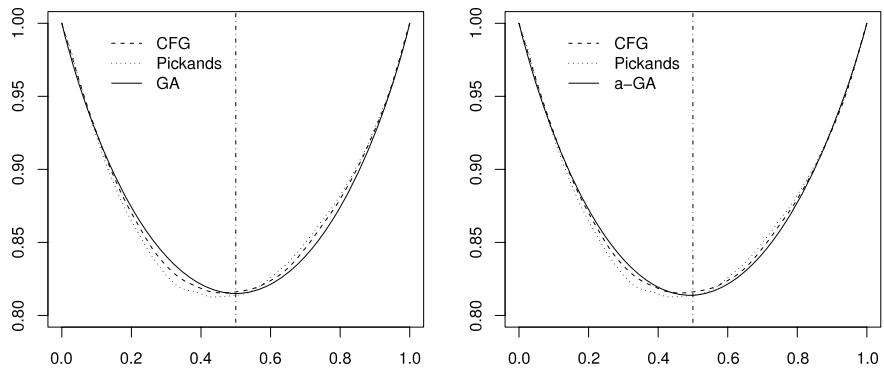}

\caption{Nonparametric estimates $A_{n,c}^{\mathrm{P}}$ and $A_{n,c}^{\mathrm
{CFG}}$, and fitted Pickands
dependence function, for the Galambos copula (left) and asymmetric
Galambos copula (right).}
\label{fig:5}
\end{figure}

Figure \ref{fig:5} displays the end-point-corrected estimates
$A_{n,c}^{\mathrm{P}}$
and $A_{n,c}^{\mathrm{CFG}}$ for the data at hand. For comparison, the
best-fitting
symmetric and asymmetric Galambos extreme-value copulas are
superimposed. Although these two models yield the highest $P$-values,
they are not significantly better than the alternatives listed in
Table \ref{tab:5}. Given the estimators' sampling variability, the data
set is simply too small to distinguish between them. This is not a
major concern, however, as predictions derived from these various
models would be roughly the same. To paraphrase Box and Draper (\cite{BoxDraper1987}, page 424), it may be that all these models are false,
but they are nearly equivalent and probably equally useful.

%
\begin{table}
\caption{Values of the statistics $S_n^{\mathrm{P}}$, $S_n^{\mathrm
{CFG}}$ and approximate
$P$-values computed using $N=2500$ parametric bootstrap samples for the
insurance data}
\label{tab:5}
\begin{tabular*}{250pt}{@{\extracolsep{\fill}}lcccc@{}}
\hline
Model & $S_n^{\mathrm{P}}$ & $P$-value & $S_n^{\mathrm{CFG}}$ &
$P$-value \\
\hline
GH & 0.087 & 0.073 & 0.048 & 0.171 \\
GA & 0.084 & 0.074 & 0.045 & 0.184 \\
HR & 0.088 & 0.067 & 0.049 & 0.157 \\
t-EV & 0.088 & 0.069 & 0.048 & 0.166 \\
a-GH & 0.052 & 0.274 & 0.012 & 0.152 \\
a-GA & 0.046 & 0.325 & 0.009 & 0.244 \\
a-HR & 0.051 & 0.272 & 0.011 & 0.174 \\
a-t-EV & 0.062 & 0.204 & 0.015 & 0.122 \\
\hline
\end{tabular*}
\end{table}

\begin{appendix}

\section*{\texorpdfstring{Appendix A: Proof of
Proposition~\protect\ref{prop:1}}{Appendix A: Proof of Proposition 1}}
\label{appa}

Let $\mathbb{A}_n$ denote either $\mathbb{A}_n^{\mathrm{P}}$ or
$\mathbb{A}_n^{\mathrm{CFG}}$ and write $\mathbb{A}_{n,\theta
_n}=\mathbb{A}_n - \mathbb{B}_{n,\theta_0}$, where $\mathbb
{B}_{n,\theta_0} = \sqrt{n}
(A_{\theta_n} - A_{\theta_0})$. As the sequence $\Theta_n$ is assumed
to converge weakly, it is tight. Thus, for given $\delta> 0$, there
exists $L = L(\delta)$ such that $\Pr( \|\Theta_n \| > L
)
$ $< \delta$ holds for every integer $n$. Therefore, for given $\zeta
> 0$,
\begin{eqnarray*}
&&\Pr\Bigl\{ \sup_{t \in[0,1]} | \mathbb{B}_{n,\theta_0} (t) - \dot
A_{\theta
_0}^\top(t) \Theta_n | > \zeta\Bigr\}\\
 &&\quad\le\Pr\Bigl\{ \sup_{t \in[0,1]} | \mathbb{B}_{n,\theta_0} (t) -
{\dot A}_{\theta
_0}^\top(t) \Theta_n | > \zeta, \| \Theta_n \| \le L \Bigr\} +
\Pr
( \|\Theta_n \| > L ) \\
&&\quad\le\Pr\Bigl\{ \sup_{t \in[0,1]} | \mathbb{B}_{n,\theta_0} (t) -
{\dot A}_{\theta
_0}^\top(t) \Theta_n | > \zeta, \| \Theta_n \| \le L \Bigr\} +
\delta.
\end{eqnarray*}

An application of the mean value theorem then implies that for every
realization $\omega$ of the process and every $t\in[0,1]$, $\mathbb
{B}_{n, \theta_0}(t,\omega) = {\dot A}^\top_{\Theta_n^*(t,\omega)}
(t)\Theta_n(\omega)$, where $\Theta_n^*(t,\omega) = \theta_0 +
\epsilon
(t,\omega) n^{-1/2} \Theta_n(\omega)$ for some $\epsilon(t,\omega)
\in
[0,1]$. It then follows from condition \eqref{eq:6} that
\begin{eqnarray*}
&&\lim_{n \to\infty} \Pr\Bigl\{ \sup_{t \in[0,1]} | \mathbb{B}_{n,
\theta_0} (t) - {\dot A}_{\theta_0}^\top(t) \Theta_n | > \zeta, \|
\Theta_n \|\le L \Bigr\}\\
&&\quad \le\lim_{n \to\infty} \Pr\Bigl\{ \| \Theta_n\| \sup_{t \in[0,1]} \big|
{\dot A}_{\Theta_n^*(t)} (t) - {\dot A}_{\theta_0} (t) | > \zeta,\|
\Theta_n \| \le L \Bigr\} \\
&&\quad\le\lim_{n \to\infty} \Pr\Bigl\{ \sup_{\|\theta- \theta_0\| \le
n^{-1/2} L} \sup_{t \in[0,1]} |{\dot A}_{\theta} (t) - {\dot
A}_{\theta
_0} (t) | > \zeta/L \Bigr\} = 0.
\end{eqnarray*}
This completes the argument.

\section*{Appendix B: Validity of the parametric bootstrap}
\label{appb}

To avoid repetitions, let $\mathbb{A}_n$ denote either $\mathbb
{A}_n^{\mathrm{P}}$ or $\mathbb{A}_n^{\mathrm{CFG}}$ and
let $\mathbb{A}$ stand for either $\mathbb{A}^{\mathrm{P}}$ or
$\mathbb{A}^{\mathrm{CFG}}$.
The following conditions, adapted from \cite{GenRem08}, ensure the
validity of the parametric bootstrap for computing $P$-values for the
proposed tests.

\begin{enumerate}[(a)]
\item[(a)]
The family $\{C_\theta\dvtx \theta\in\mathcal{O} \}$ of extreme-value
copulas must be such that:
\begin{enumerate}[(viii)]
\item[(i)]
the parameter space $\mathcal{O} $ is an open subset of $\mathbb{R}^p$;
\item[(ii)]
members of the family are identifiable, that is, for
every $\epsilon> 0$,
\[
\inf\Bigl\{ \sup_{t \in[0,1]} \| A_\theta(t) - A_{\theta_0} (t)
\| \dvtx
\theta\in\mathcal{O} \mbox{ and } \|\theta- \theta_0\| > \epsilon
\Bigr\} > 0;
\]

\item[(iii)] the mapping $\theta\mapsto A_\theta$ is Fr\'echet
differentiable with derivative $\theta\mapsto\dot A_\theta$, that is,
for all $\theta_0 \in\mathcal{O} $,
\[
\lim_{\|h\| \downarrow0} \sup_{t \in[0,1]} \frac{\|A_{\theta_0 + h}(t)
- A_{\theta_0} (t) - \dot A_{\theta_0}^\top(t)h \|}{\| h \| } = 0;
\]

\item[(iv)] $C_\theta$ has a Lebesgue density $c_\theta$ for
all $\theta\in\mathcal{O} $;

\item[(v)] the density $c_\theta$ admits first- and second-order derivatives
with respect to all components of $\theta\in\mathcal{O} $; the
gradient (row) vector with respect to $\theta$ is denoted $\dot
c_\theta
$ and the Hessian matrix is denoted $\ddot c_\theta$;

\item[(vi)] for arbitrary $(u,v) \in(0,1)^2$ and every $\theta_0 \in
\mathcal{O} $, $ \theta\mapsto\dot
c_\theta(u,v)/c_{\theta}(u,v)$ and $\theta\mapsto\ddot
c_\theta(u,v)/c_{\theta}(u,v)$ are continuous at $\theta_0$,
$C_{\theta_0}$ almost surely;

\item[(vii)] for every $\theta_0 \in\mathcal{O} $, there exist a
neighborhood
$\mathcal{N}$ of $\theta_0$ and a Lebesgue integrable function $h\dvtx
(0,1)^2 \to\mathbb{R}$ such that $
\sup_{\theta\in\mathcal{N}} \| \dot c_\theta(u,v) \|
\le h (u,v)$ holds for all $(u,v) \in(0,1)^2$;

\item[(viii)] for every $\theta_0 \in\mathcal{O} $, there exist a
neighborhood
$\mathcal{N}$ of $\theta_0$ and $C_{\theta_0}$-integrable functions
$h_1, h_2\dvtx  (0,1)^2 \to\mathbb{R}$ such that for all $(u,v) \in
(0,1)^2$,
\[
\sup_{\theta\in\mathcal{N}} \biggl\| \frac{\dot
c_{\theta}(u,v)}{c_{\theta}(u,v)}\biggr \|^2 \le h_1(u,v)
\quad \mbox{and}\quad  \sup_{\theta\in\mathcal{N}} \biggl\|\frac{\ddot
c_{\theta}(u,v)}{c_{\theta}(u,v)} \biggr\| \le h_2(u,v).
\]
\end{enumerate}

\item[(b)] In addition, the estimators $A_n$ and $\theta_n$ satisfy
the following:

\begin{enumerate}[(ii)]

\item[(i)]$(\mathbb{A}_n,\Theta_n,\mathbb{W}_n)\rightsquigarrow
(\mathbb{A},\Theta
,\mathbb{W})$ in $\mathcal{D}([0,1],\mathbb{R})\times\mathbb
{R}^{p\otimes2}$ as $n \to\infty$, where the limit is a centered
Gaussian process. Here,
\[
\mathbb{W}_n = n^{-1/2} \sum_{i=1}^n \frac{\dot c_{\theta_0}^\top
(U_i^*, V_i^*)}{c_{\theta_0} (U_i^*, V_i^*)}
\]
for a random sample $(U_1^*,V_1^*),\dots, (U_n^*,V_n^*)$ from
$C_{\theta
_0}$ and $\mathbb{W}$ is $\mathcal{N}(0,I_P)$, where $I_P$ is the
Fisher information matrix; see \cite{GenRem08}, page 1101.

\item[(ii)] $\mathrm{E}_{\theta_0} (\Theta\mathbb{W}^\top)= J$, where $J$
is the $p\times p$ identity matrix. Further, $\mathrm{E}_{\theta_0} \{
\mathbb{A}
(t)\mathbb{W}\} = \dot A_{\theta_0}(t)$ for every $t\in(0,1)$.

\end{enumerate}

\end{enumerate}

Condition (b) can be checked as follows, under the assumption that
$(\mathbb{A}
_n,\Theta_n) \rightsquigarrow(\mathbb{A},\Theta)$ as $n \to\infty
$. First,
results from Chapter 5 of \cite{GanStu87} can be combined with the
functional delta method (see, e.g., \cite{vanWel96}, Section 3.9) to
see that as $n \to\infty$, $(\mathbb{A}_n,\Theta_n,\mathbb
{C}_n,\mathbb{W}_n)
\rightsquigarrow(\mathbb{A},\Theta,\mathbb{C},\mathbb{W})$.

Next, observe that $ \mathrm{E}_{\theta_0} \{\mathbb{C}(u,v)\mathbb
{W}\} = \dot
C_{\theta_0}(u,v)$ for all $u, v \in[0,1]$; see \cite{GenRem08}, page 1108. Given that, for all $t \in[0,1]$,
\[
\mathbb{A}^{\mathrm{P}}(t) = -A^2_{\theta_0}(t) \int_0^1 \mathbb
{C}_{\theta_0}(x^{1-t}, x^t)
\frac{ \mathrm{d}x}{x}
\]
and
\[
\mathbb{A}^{\mathrm{CFG}}(t) = A_{\theta_0}(t) \int_0^1 \mathbb
{C}_{\theta_0}(x^{1-t}, x^t)
\frac{ \mathrm{d}x}{x\log(x)} ,
\]
we can see that
\[
\mathrm{E}_{\theta_0} \{ \mathbb{A}^\mathrm{P}(t)\mathbb{W} \} =
-A_{\theta
_0}^2(t)\int_0^1 \dot C_{\theta_0}(x^{1-t}, x^t) \frac{\mathrm{d}x}{x}
\]
and
\[
\mathrm{E}_{\theta_0} \{ \mathbb{A}^\mathrm{CFG}(t)\mathbb{W} \} =
A_{\theta
_0}(t)\int_0^1 \dot C_{\theta_0}(x^{1-t}, x^t) \frac{ \mathrm
{d}x}{x\log
(x)} .
\]

Interchanging the order of differentiation and integration, we get
$\mathrm{E}_{\theta_0} \{ \mathbb{A}^\mathrm{P}(t)\mathbb{W} \} =\break
\mathrm{E}_{\theta_0}
\{ \mathbb{A}^\mathrm{CFG}(t)\mathbb{W} \} = \dot A_{\theta_0}(t)$ for
all $t \in(0,1)$.

As for the condition $\mathrm{E}_{\theta_0}(\Theta\mathbb{W}) = J$, it can
be verified using \cite{GenRem08}, Proposition 4, for the estimators
based on maximum pseudo-likelihood and on the inversion of Spearman's
rho. To handle the estimator based on Kendall's tau,
Proposition 5 in \cite{GenRem08} must be used instead.

\section*{\texorpdfstring{Appendix C: Proof of Proposition
\protect\ref{prop:2}}{Appendix C: Proof of Proposition 2}}
\label{appc}

The proof closely mimics the argument presented in \cite{GenSeg09},
Appendix B. To avoid duplication, the same notation is used and only
the critical differences are highlighted. This also offers an
opportunity to correct minor typographical errors in the original source.

First, consider the process given by $\mathbb{B}_{n}^{\mathrm{P}}(t)
= n^{1/2} \{1/A_n^{\mathrm{P}}
(t) - 1/A_C^\mathrm{P}(t) \}$ for all $t \in[0,1]$ and show that
$\mathbb{B}_{n}^{\mathrm{P}}\rightsquigarrow\mathbb{B} = - \mathbb
{A}^{\mathrm{P}}_C /(A_C^\mathrm{P})^2$ as $n
\to\infty$. Then
\[
\sqrt{n} (A_n^{\mathrm{P}}- A_C^\mathrm{P}) = \frac{-(A_C^\mathrm{P})^2
\mathbb{B}_{n}^{\mathrm{P}}
}{1+n^{-1/2}\mathbb{B}_{n}^{\mathrm{P}}A^\mathrm{P}_C} \rightsquigarrow
\mathbb{A}^{\mathrm{P}}_C,
\]
as a consequence of the functional version of Slutsky's lemma.

Put $k_n = 2 \log(n+1)$ and write
\[
\mathbb{B}_{n}^{\mathrm{P}}(t) = \int_0^1 \mathbb{C}_n(x^{1-t},
x^{t}) \frac{\mathrm{d}x}{x}
=\int_0^\infty\mathbb{C}_n\bigl(\mathrm{e}^{-s(1-t)}, \mathrm{e}^{-st}\bigr)\, \mathrm{d}s =
I_{1,n} + I_{2,n},
\]
where, for each $t \in[0,1]$,
\[
I_{1,n} (t) = \int_{k_n}^\infty\mathbb{C}_n\bigl(\mathrm{e}^{-s(1-t)}, \mathrm{e}^{-st}\bigr)
\,\mathrm{d}s ,\qquad I_{2,n} (t) = \int_0^{k_n} \mathbb{C}_n\bigl(\mathrm{e}^{-s(1-t)},
\mathrm{e}^{-st}\bigr)\, \mathrm{d}s.
\]
The contribution of $I_{1,n}(t)$ is asymptotically negligible because
the fact that $s > k_n$ implies that $\min(\mathrm{e}^{-s(1-t)} , \mathrm{e}^{-st}
) < 1/(n+1)$ and hence
that
\[
\bigl|\mathbb{C}_n\bigl(\mathrm{e}^{-s(1-t)}, \mathrm{e}^{-st}\bigr)\bigr| = n^{1/2} C\bigl(\mathrm{e}^{-s(1-t)}, \mathrm{e}^{-st}\bigr)
\le n^{1/2} \min\bigl(\mathrm{e}^{-s(1-t)}, \mathrm{e}^{-st}\bigr) \le n^{1/2} \mathrm{e}^{-s/2}.
\]
Thus, for all $t \in[0,1]$,
%
%
\setcounter{equation}{0}
\begin{equation}
\label{eq:A1}
| I_{1,n}(t)| \le n^{1/2} \int_{k_n}^\infty C\bigl(\mathrm{e}^{-s(1-t)}, \mathrm{e}^{-st}\bigr)
\,\mathrm{d}s \\
\leq n^{1/2} \int_{k_n}^\infty \mathrm{e}^{-s/2} \,\mathrm{d}s \leq\frac
{2}{n^{1/2}}.
\end{equation}

Consequently, the asymptotic behavior of $\mathbb{B}_{n}^{\mathrm
{P}}$ is determined
entirely by $I_{2,n}$. Invoking the Stute representation given by
Genest and Segers \cite
{GenSeg09}, we may write $I_{2,n} = J_{1,n} + J_{2,n} + J_{3,n} +
\mathrm{o}(1)$, where, for each $t \in[0,1 ]$,
\begin{eqnarray*}
J_{1,n} (t) &=& \int_0^{k_n} \alpha_n\bigl(\mathrm{e}^{-s(1-t)}, \mathrm{e}^{-st}\bigr)
\,\mathrm{d}s,
\\
J_{2,n} (t) &=& - \int_0^{k_n} \alpha_n\bigl(\mathrm{e}^{-s(1-t)}, 1\bigr)
{\dot C}_1\bigl(\mathrm{e}^{-s(1-t)}, \mathrm{e}^{-st}\bigr)\, \mathrm{d}s , \\
J_{3,n} (t) &=& - \int_0^{k_n} \alpha_n(1, \mathrm{e}^{-st}) {\dot
C}_2\bigl(\mathrm{e}^{-s(1-t)}, \mathrm{e}^{-st}\bigr)\, \mathrm{d}s .
\end{eqnarray*}
Here, ${\dot C}_1 (u,v) = \partial C(u,v)/\partial u$, ${\dot C}_2
(u,v) = \partial C(u,v)/\partial v$ and $\alpha_n$ is the empirical
process associated with the pairs $(F(X_1),G(Y_1)), \ldots, (F(X_n),G(Y_n))$.

Fix $\omega\in(0, 1/2)$ and write $q_\omega(t) = t^\omega
(1-t)^\omega$ for all $t \in[0,1]$. Also, let
\begin{eqnarray*}
K_1 (s,t) &=& q_\omega\bigl\{ \min\bigl( \mathrm{e}^{-s(1-t)} , \mathrm{e}^{-st}\bigr) \bigr\}, \\
K_2 (s,t) &=& q_\omega\bigl(\mathrm{e}^{-s(1-t)}\bigr) {\dot C}_1\bigl(\mathrm{e}^{-s(1-t)},
\mathrm{e}^{-st}\bigr), \\
K_3 (s,t) &=& q_\omega(\mathrm{e}^{-st}) {\dot C}_2\bigl(\mathrm{e}^{-s(1-t)}, \mathrm{e}^{-st}\bigr)
\end{eqnarray*}
for all $s \in(0, \infty)$ and $t \in[0,1]$. The proof that $J_{1,n}
+ J_{2,n} + J_{3,n}$ has the stated limit then proceeds exactly as in
Appendix B of \cite{GenSeg09}, provided that for $i=1,2,3$, there
exists an integrable function $K_i^*\dvtx  (0,\infty) \to\mathbb{R}$ such
that $K_i (s,t) \le K_i^* (s)$ for all $s \in(0, \infty)$ and $t \in[0,1]$.

For $K_1$, this is immediate because $K_1 (s,t) \leq \mathrm{e}^{-\omega s/2}$
for all $s \in(0, \infty)$ and $t \in[0,1]$. For $K_2$, the facts
that $C$ is LTD and smaller than the Fr\'echet--Hoeffding upper bound
imply that
\[
{\dot C}_1 \bigl(\mathrm{e}^{-s(1-t)}, \mathrm{e}^{-st}\bigr) \le \mathrm{e}^{s(1-t)} C \bigl(\mathrm{e}^{-s(1-t)},
\mathrm{e}^{-st}\bigr) \le \mathrm{e}^{s(1-t)} \min\bigl(\mathrm{e}^{-s(1-t)}, \mathrm{e}^{-st}\bigr).
\]
Now, set $m(t) = \max(t,1-t)$ and note that $q_\omega(\mathrm{e}^{-s(1-t)})
\le
\mathrm{e}^{-\omega s(1-t)}$ for all $s \in(0, \infty)$ and $t \in[0,1]$. Therefore,
\[
K_2 (s,t) \leq \mathrm{e}^{s(1-\omega)(1-t)} \mathrm{e}^{-sm(t)} \leq \mathrm{e}^{s(1-\omega)m(t)}
\mathrm{e}^{-sm(t)}
= \mathrm{e}^{-s\omega m(t) } \le \mathrm{e}^{-\omega s/2}
\]
because $m(t) \ge1/2$ for all $t \in[0,1]$. The argument for $K_3$ is similar.

Turning to the $A_n^{\mathrm{CFG}}$ estimator, observe that
\begin{eqnarray*}
\mathbb{B}_{n}^{\mathrm{CFG}}(t) &=& n^{1/2} \{ \log A_{n}^{\mathrm
{CFG}}(t) - \log A^\mathrm{CFG}_C (t) \}
\\ &=&
\int_0^1 \mathbb{C}_n(x^{1-t}, x^t) \frac{\mathrm{d}x}{x \log(x)}
= -
\int_0^\infty\mathbb{C}_n\bigl(\mathrm{e}^{-s(1-t)}, \mathrm{e}^{-st}\bigr) \frac{\mathrm{d}s}{s}
\end{eqnarray*}
for all $t \in[0,1]$. This process can be written as
$-(I_{1,n}+I_{2,n}+I_{3,n})$,
where
\begin{eqnarray*}
I_{1,n} (t) &=& \int_{k_n}^\infty\mathbb{C}_n\bigl(\mathrm{e}^{-s(1-t)},
\mathrm{e}^{-st}\bigr) \frac{\mathrm{d}s}{s} , \\
I_{2,n} (t) &=& \int_{\ell_n}^{k_n} \mathbb{C}_n\bigl(\mathrm{e}^{-s(1-t)},
\mathrm{e}^{-st}\bigr) \frac{\mathrm{d}s}{s} , \\
I_{3,n} (t) &=& \int_0^{\ell_n} \mathbb{C}_n\bigl(\mathrm{e}^{-s(1-t)}, \mathrm{e}^{-st}\bigr)
\frac{\mathrm{d}s}{s}
\end{eqnarray*}
with $k_n = 2 \log(n+1)$ as above and $\ell_n = 1/(n+1)$.

Arguing as in \eqref{eq:A1}, we see that $|I_{1,n}| \le n^{-1/2}$.
Similarly, $I_{3,n}$ is negligible asymptotically, for if $s \in(0,
\ell_n)$ and $t \in[0, 1]$, then we have
\[
\min\bigl(\mathrm{e}^{-s(1-t)}, \mathrm{e}^{-st}\bigr) \geq \mathrm{e}^{-1/(n+1)} > \frac{n}{n+1}
\]
and hence $C_n (\mathrm{e}^{-s(1-t)}, \mathrm{e}^{-st}) = 1$. Furthermore, the fact that
$C$ is LTD implies that
$C (\mathrm{e}^{-s(1-t)}, \mathrm{e}^{-st}) \ge \mathrm{e}^{-s}$ for all $s \in(0, \infty)$ and
$t \in[0,1]$. Therefore,
\[
\bigl|\mathbb{C}_n\bigl(\mathrm{e}^{-s(1-t)}, \mathrm{e}^{-st}\bigr) \bigr| \leq n^{1/2} ( 1 - \mathrm{e}^{-s} ) \leq
n^{1/2} s.
\]
Consequently, $|I_{3,n}| \le n^{1/2} \ell_n \leq n^{-1/2}$. As a result,
the asymptotic behavior of $\mathbb{B}_{n}^{\mathrm{CFG}}$ is
determined entirely by
$I_{2,n}$. Following Genest and Segers \cite{GenSeg09}, we can further write $I_{2,n} =
J_{1,n} + J_{2,n} + J_{3,n} + \mathrm{o}(1)$, where, for all $t \in[0,1 ]$,
\begin{eqnarray*}
J_{1,n} (t) &=& \int_{\ell_n}^{k_n} \alpha_n\bigl(\mathrm{e}^{-s(1-t)}, \mathrm{e}^{-st}\bigr)
\frac{\mathrm{d}s}{s} , \\
J_{2,n} (t) &=& - \int_{\ell_n}^{k_n} \alpha_n\bigl(\mathrm{e}^{-s(1-t)},1\bigr) {\dot
C}_1\bigl(\mathrm{e}^{-s(1-t)}, \mathrm{e}^{-st}\bigr) \frac{\mathrm{d}s}{s} , \\
J_{3,n} (t) &=& - \int_{\ell_n}^{k_n} \alpha_n(1, \mathrm{e}^{-st}) {\dot
C}_2\bigl(\mathrm{e}^{-s(1-t)}, \mathrm{e}^{-st}\bigr) \frac{\mathrm{d}s}{s} .
\end{eqnarray*}

The joint asymptotic behavior of these terms can be determined in the
same way as before. The only difference is that the integration measure
is now $\mathrm{d}s/s$. For $s \in[1, \infty)$, the same upper bounds
$K_1^*$, $K_2^*$, $K_3^*$ apply and they have already been shown to be
integrable on this domain. To obtain an integrable bound for $K_1$ on
$(0,1)$, it suffices to use the fact that $K_1 (s,t) \le(1-
\mathrm{e}^{-sm(t)})^\omega\le\{ sm(t) \}^\omega\le s^\omega$. The same bound
works for both $K_2$ and $K_3$ because $\dot C_i \in[0,1]$ for
$i=1,2$. This completes the argument. 
\end{appendix}

\section*{Acknowledgements}  

This research was supported by grants from the Natural Sciences and
Engineering Research Council of Canada, the Fonds qu\'eb\'ecois de
la recherche sur la nature et les technologies, and the
Institut de finance math\'ematique de montr\'eal.
Some of the computations were carried out on a Beowulf cluster at the
Department of Statistics, University of Connecticut, which was partially
supported by a grant from the National Science Foundation, Scientific
Computing Research Environments for the Mathematical Sciences (SCREMS) Program.  

\printhistory

\end{document}